\documentclass[]{elsarticle} 
\usepackage[hyphens]{url}

\usepackage{lineno} 
\providecommand{\tightlist}{%
  \setlength{\itemsep}{0pt}\setlength{\parskip}{0pt}}

\usepackage{graphicx}
\usepackage{booktabs} 

\usepackage[T1]{fontenc}
\usepackage{lmodern}
\usepackage{amssymb,amsmath}
\usepackage{ifxetex,ifluatex}
\usepackage{fixltx2e} 
\IfFileExists{upquote.sty}{\usepackage{upquote}}{}
\ifnum 0\ifxetex 1\fi\ifluatex 1\fi=0 
  \usepackage[utf8]{inputenc}
\else 
  \usepackage{fontspec}
  \ifxetex
    \usepackage{xltxtra,xunicode}
  \fi
  \defaultfontfeatures{Mapping=tex-text,Scale=MatchLowercase}
  
\fi
\IfFileExists{microtype.sty}{\usepackage{microtype}}{}
\ifxetex
  \usepackage[setpagesize=false, 
              unicode=false, 
              xetex]{hyperref}
\else
  \usepackage[unicode=true]{hyperref}
\fi
\hypersetup{breaklinks=true,
            bookmarks=true,
            pdfauthor={},
            pdftitle={Statistical Inference for inter-arrival times of extreme events in bursty time series},
            colorlinks=false,
            urlcolor=blue,
            linkcolor=magenta,
            pdfborder={0 0 0}}
\begin{document}
\begin{frontmatter}

  \title{Statistical Inference for inter-arrival times of extreme events in
bursty time series}
    \author[a]{Katharina Hees}
   \ead{hees@statistik.tu-dortmund.de} 
    \author[b]{Smarak Nayak}
   \ead{smarak.nayak@nab.com.au} 
    \author[c]{Peter Straka}
   \ead{straka.ps@gmail.com} 
      \address[a]{Department of Statistics, TU Dortmund University, Dortmund, Germany}
    \address[b]{National Australia Bank, Melbourne, Australia}
    \address[c]{School of Mathematics and Statistics, UNSW, Sydney, Australia}
    
  \begin{abstract}
  In many complex systems studied in statistical physics, inter-arrival
  times between events such as solar flares, trades and neuron voltages
  follow a heavy-tailed distribution. The set of event times is
  fractal-like, being dense in some time windows and empty in others, a
  phenomenon which has been dubbed ``bursty''. A new model for the
  \emph{inter-exceedance times} of such events above high thresholds is
  proposed. For high thresholds and infinite-mean waiting times, it is
  shown that the times between threshold crossings are Mittag-Leffler
  distributed, and thus form a ``fractional Poisson Process'' which
  generalizes the standard Poisson Process of threshold exceedances.
  Graphical means of estimating model parameters and assessing model fit
  are provided. The inference method is applied to an empirical bursty
  time series, and it is shown how the memory of the Mittag-Leffler
  distribution affects prediction of the time until the next extreme
  event."
  \end{abstract}
 \end{frontmatter}

\hypertarget{introduction}{%
\section{Introduction}\label{introduction}}

Time series displaying temporally inhomogeneous behaviour in terms of
the occurrence of events have received strong interest in the recent
statistical physics literature (Bagrow and Brockmann, 2013; Barabási,
2005; Karsai et al., 2011; Min et al., 2011; Oliveira and Barabási,
2005; Omi and Shinomoto, 2011; Vasquez et al., 2006; Vazquez et al.,
2007). They have been observed in the context of earthquakes, sunspots,
neuronal activity and human communication (see Karsai et al., 2012;
Meerschaert and Stoev, 2008 for a list of references; Vajna et al.,
2013). Such time series exhibit high activity in some `bursty'
intervals, which alternate with other, quiet intervals. Although several
mechanisms are plausible explanations for bursty behaviour (most
prominently self-exciting point process by Hawkes (1971) and renewal
Hawkes processes, e.g.~Wheatley et al. (2016), Stindl and Chen (2018)),
there seems to be one salient feature which very typically indicates the
departure from temporal homogeneity: a heavy-tailed distribution of
waiting times (Karsai et al., 2012; Vajna et al., 2013; Vasquez et al.,
2006). As we show below in simulations, a simple renewal process with
heavy-tailed waiting times can capture this type of dynamics. For many
systems, the renewal property is appropriate; a simple test of the
absence of correlations in a succession of waiting times can be
undertaken by randomly reshuffling the waiting times (Karsai et al.,
2012).

Often a magnitude, or mark can be assigned to each event in the renewal
process, such as for earthquakes, solar flares or neuron voltages. The
Peaks-Over-Threshold model (POT, see e.g.~Coles, 2001) applies a
threshold to the magnitudes, and fits a Generalized Pareto distribution
to the threshold exceedances. A commonly made assumption in POT models
is that times between events are either fixed or light-tailed, and this
entails that the threshold crossing times form a Poisson process (Hsing
et al., 1988). Then as one increases the threshold \(\ell\) and thus
decreases the threshold crossing probability \(p_{\ell}\), the Poisson
process is thinned, i.e.~its intensity decreases \emph{linearly} with
\(p_{\ell}\) (see e.g.~Beirlant et al., 2006).

As will be shown below, in the heavy-tailed waiting time scenario
threshold crossing times form a \emph{fractional Poisson process}
(Laskin, 2003; Meerschaert et al., 2011), which is a renewal process
with Mittag-Leffler distributed waiting times. The family of
Mittag-Leffler distributions nests the exponential distribution (Haubold
et al., 2011), and hence the fractional Poisson process generalizes the
standard Poisson process. Again as the threshold size \(\ell\) increases
and the threshold crossing probability \(p_{\ell}\) decreases, the
fractional Poisson process is thinned: The scale parameter of the
Mittag-Leffler inter-arrival times of threshold crossing times
increases, but \emph{superlinearly}; see the Theorem below.

Maxima of events which occur according to a renewal process with
heavy-tailed waiting times have been studied under the names
``Continuous Time Random Maxima process'' (CTRM) (Benson et al., 2007;
Hees and Scheffler, 2018a, 2018b; Meerschaert and Stoev, 2008),
``Max-Renewal process'' (Basrak and Špoljarić, 2015; Silvestrov, 2002;
Silvestrov and Teugels, 2004), and ``Shock process'' (Anderson, 1987;
Esary and Marshall, 1973; Gut and Hüsler, 1999; Shanthikumar and Sumita,
1985, 1984, 1983). The existing literature focuses on probabilistic
results surrounding these models. In this work, however, we introduce a
method of inference for this type of model, which is seemingly not
available in the literature.

We review the marked renewal process in Section 2.1, and derive a
scaling limit theorem for inter-exceedance times in Section 2.2. We give
a statistical procedure to estimate model parameters via stability plots
in Section 3.1 and 3.4, but first we need to discuss inference for the
Mittag-Leffler distribution in Section 3.2 as well as a Likelihood-ratio
test to guide the choice, whether a Mittag-Leffler or an exponential
distribution fits better to the inter-exceedance times. A simulation
study of the effectiveness of our statistical procedure is given in
Section 4. In Section 5 we apply our method to a real data set. In
Section 6, we discuss the memory property of the Mittag-Leffler
distribution, and how it affects the predictive distribution for the
time until the next threshold crossing event. Finally we close with a
discussion and conclusion in Section 7. For all statistical computations
we have used \texttt{R} (R Core Team, 2018) and the package
\texttt{CTRE} (Hees and Straka, 2018). Source code for simulations and
figures generated in this manuscript is available online at
\url{https://github.com/strakaps/bursty-POT}.

\hypertarget{probabilistic-model}{%
\section{Probabilistic Model}\label{probabilistic-model}}

\hypertarget{sec:CTRE}{%
\subsection{Continuous Time Random Exceedances (CTRE)}\label{sec:CTRE}}

A Borel-measurable function \(f:\mathbb (0,\infty) \to (0,\infty)\) is
said be ``regularly varying'' at \(\infty\) with parameter (or
``index'') \(\rho\) if

\[
\lim_{x \to \infty} \frac{f(\lambda x)}{f(x)} = \lambda^\rho \quad \text{ for all } \lambda > 0.
\]

For more details on regular variation and stable limit theorems, we
recommend the book by Meerschaert and Sikorskii (2012).

As a model for extreme observations, we use a Marked Renewal Process
(MRP):

\begin{description}
\item[\textbf{Definition (MRP):}]
Let \((W,J), (W_1, J_1), (W_2, J_2), \ldots\) be iid pairs of random
variables, where the \(W_k > 0\) are interpreted as the \emph{waiting
times} and \(J_k \in \mathbb{R}\) as the \emph{event magnitudes}. If
\(W\) and \(J\) are independent, the Marked Renewal Process is said to
be \emph{uncoupled}. We call the MRP \emph{bursty}, if the tail function
\(F_W(x)=P(W>t)\) is regularly varying with index \(0<\beta<1\).
\end{description}

In the following we denote with \(x_L \in [-\infty, +\infty)\) and
\(x_R \in (-\infty, +\infty]\) the left and right endpoint of the
distribution of \(J\). We assume that the \(k\)-th magnitude \(J_k\)
occurs at time \(T_k = W_1 + \ldots + W_k\). Based on a MRP, we define
the Continuous Time Random Exceedance model (CTRE) as follows:

\begin{description}
\item[\textbf{Definition (CTRE):}]
Given a threshold \(\ell \in (x_L, x_R)\), consider the stopping time
\[\tau(\ell) := \min\{k: J_k > \ell\},\quad \ell \in (x_L, x_R).\]
Define the pair of random variables \((X(\ell), T(\ell))\) via
\[X(\ell) = J_{\tau(\ell)} - \ell, \quad 
T(\ell) = \sum_{k=1}^{\tau(\ell)} W_k.\] By restarting the MRP at
\(\tau(\ell)\), inductively define the two iid sequences \(T(\ell,n)\)
and \(X(\ell, n)\), \(n \in \mathbb N\), called the ``inter-exceedance
times'' (IETs) and the ``exceedances'', respectively. The pair sequence
\((X(\ell, n), T(\ell, n))_{n \in \mathbb N}\) is called a Continuous
Time Random Exceedance model (CTRE). If the underlying MRP is uncoupled,
then the CTRE is also called uncoupled. We call the CTRE \emph{bursty},
if the tail function of the IETs is regularly varying with an index
\(0<\beta<1\).
\end{description}

\begin{figure}

{\centering \includegraphics[width=0.7\linewidth]{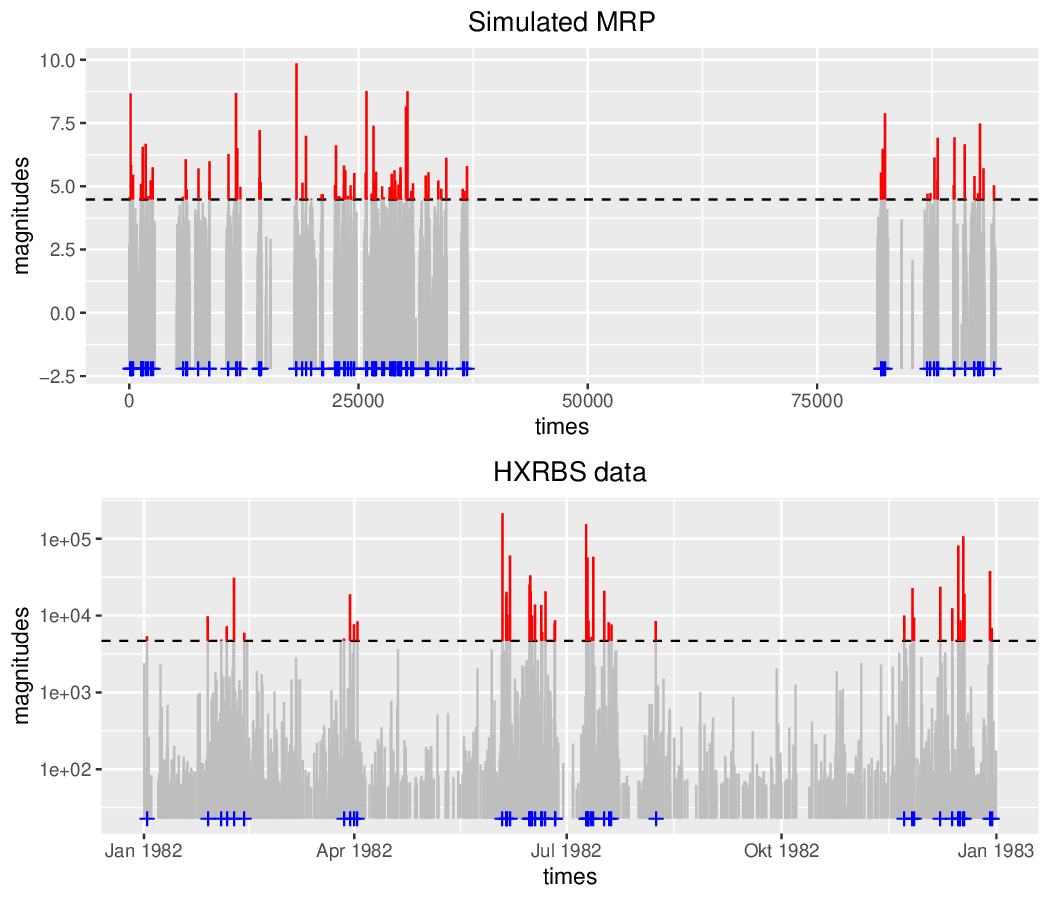} 

}

\caption{\label{fig:thresholdedBursty} Exceedances (red) and times until Exceedance (durations between blue crosses) for a given threshold $\ell$ (dashed line). Upper picture: Simulated data with stable distributed waiting times. Lower picture: Solar flares during 1982.}\label{fig:thresholdedBursty}
\end{figure}

In this article, we restrict ourselves to the uncoupled case, where
\(W\) and \(J\) are independent. Then the two sequences
\(X(\ell, n)_{n \in \mathbb N}\) and \(T(\ell, n)_{n \in \mathbb N}\)
are independent as well. To see why, note that \(X(\ell)\) is, in
distribution, simply equal to \(J - \ell | J > \ell\), independent of
any waiting time \(W_k\).\\
We assume for the rest of the article, that the magnitudes
\((J_i)_{i \in \mathbb{N}}\) belong to the the max-domain of attraction
of some non-degenerate distribution. This means there exist \(a_n>0\)
and \(d_n \in \mathbb{R}\) such that \begin{align} \label{assumptionJs}
a_n^{-1}(J_1 \vee \ldots \vee J_n -d_n) \Rightarrow A \text{ as } n \rightarrow \infty.
\end{align} Hence, the distribution of \(A\) is a generalized extreme
value distribution (GEV) whose distribution function is given by
\[ F(x;\xi) = \begin{cases}\exp(-(1+\xi x)_+^{-1/\xi}) & \xi\neq0 \\ \exp(-\exp(-x)) & \xi = 0\end{cases}\]
where \((.)_+:= \max\{0,.\}\).

The GEV is subdivided into the Gumbel (\(\xi=0\)), the Weibull
(\(\xi<0\)) and the Fréchet (\(\xi>0\)) family of distributions.

Figure \ref{fig:thresholdedBursty} shows a simulated dataset in the top
panel, where \(W\) has a stable distribution with tail parameter
\(\beta =\) 0.8 (and skewness \(1\) and location \(0\)), and where \(J\)
is from a standard Gumbel distribution. In the bottom panel, we plot a
time series of solar flare intensities derived from a NASA dataset
(Dennis et al., 1991) which we will later examine more closely (see
Section 7). Clearly, the simulated data exhibit long intervals
\emph{without any} events, whereas in the real-world dataset events
appear continuously. The threshold exceedances, however, appear to have
visually similar statistical behaviour in both models. Observations
below a threshold are commonly discarded in Extreme Value Theory (POT
approach); likewise, the CTRE model interprets these observations as
noise and discards them.

\hypertarget{sec:scaling}{%
\subsection{Scaling limit of Exceedance Times}\label{sec:scaling}}

In this section we state and prove the key theorem, which is founded on
the concept of regular variation.

\begin{description}
\item[\textbf{Theorem:}]
For the magnitudes \(J_k\), let assumption \eqref{assumptionJs} hold.
Furthermore, let the waiting times \(W_k\) be in the domain of
attraction of a positively skewed sum-stable law with stability
parameter \(0 < \beta < 1\); more precisely,
\begin{align} \label{eq:stability}
\frac{W_1 + \ldots + W_n}{b(n)} \Rightarrow D, 
\quad n \to \infty
\end{align} for a function \(b(n)\) which is regularly varying at
\(\infty\) with parameter \(1/\beta\), and where
\(\mathbf E[\exp(-sD)] = \exp(-s^\beta)\), \(s > 0\). Write
\(p_{\ell} := \mathbf P(J > \ell)\). Then the weak convergence \[
\frac{T(\ell)} {b(1/p_{\ell})} \Rightarrow  Z_\beta \quad \text{ as } \quad \ell \uparrow x_R
\] holds, where the Mittag-Leffler random variable \(Z_\beta\) is
defined on the positive real numbers via \[
\mathbf E[\exp(-sZ_\beta)] = \frac{1}{1+s^\beta}.
\]
\end{description}

\emph{Proof of Theorem:} Due to assumption \eqref{assumptionJs}
\begin{align*}
& & P\left(\frac{J_1 \vee \ldots \vee J_{\lfloor c \rfloor}-d(c)}{a(c)} \leq x\right) & \rightarrow F(x;\xi) \\
& \Longleftrightarrow & F_J(xa(c)+d(c))^{\lfloor c \rfloor} & \rightarrow F(x;\xi) \\
& \Longleftrightarrow &  \lfloor c \rfloor \log F_J(\ell_c) & \rightarrow \log F(x;\xi)
\end{align*} as \(c \rightarrow \infty\), for any \(x\) from the support
of \(F(x;\xi)\), with \(\ell_c:=xa(c)+d(c)\). Furthermore, since
\(\log(1-x) \sim -x\) for small \(x\) it follows \begin{align*}
c \cdot p(\ell_c) \rightarrow -\log F(x;\xi) \quad  \text{ as } c \rightarrow \infty,
\end{align*} with \(p(\ell_c):= 1-F(\ell_c)\). Due to
\(\tau(\ell_c) \sim {\rm Geo}(p(\ell_c))\) and above equation, it
follows that \(\tau(\ell_c)/c\) converges to an exponential random
variable, \begin{align*}
    \frac{\tau(\ell_c)}{c} \Rightarrow E \quad  \text{ as } c \rightarrow \infty
\end{align*} with inverse mean \(\lambda:=-\log F(x;\xi)\). Due to
\begin{align*}
    S(c):= \sum_{i=1}^{\lfloor c \rfloor} \frac{W_i}{b(c)} \Rightarrow D \quad \text{ as } c \rightarrow \infty,
\end{align*} it follows with Gnedenko's transfer theorem (see Gnedenko
(1983)), that \begin{align*}
    \sum_{i=1}^{ \tau(\ell_c)} \frac{W_i}{b(\ell_c)}  \Rightarrow Z \quad \text{ as }  c \rightarrow \infty, 
\end{align*} where the distribution of \(Z\) has the characteristic
function \begin{align*}
\Psi(s) &= \int_0^{\infty} (\Phi_D(s))^{y} F_E(dy) = \frac{1}{1-\log(\Phi_D(s))/\lambda}\\
&= \frac{1}{1+s^{\beta}\lambda^{-1}} = \frac{1}{1+(s \lambda^{-1/\beta})^{\beta}}
\end{align*} where \(\Phi_D(s) = \exp(-s^\beta)\) is the Laplace
transform of \(D\) and \(F_E(dy)\) is the distribution function of
\(E\). Hence, \(Z\) is Mittag-Leffler distributed with scale parameter
\(\lambda^{-1/\beta}\). Rewriting \begin{align*}
    \sum_{i=1}^{ \tau(\ell_c)} \frac{W_i}{b(1/p(\ell_c))} = \left(\sum_{i=1}^{ \tau(\ell_c)} \frac{W_i}{b(c) }\right) \frac{b(c)}{b(1/p(\ell_c))}  
\end{align*} and \begin{align*}
    \frac{b(c)}{b(1/p(\ell_c))} = \frac{b(c)}{b(c/(p(\ell_c)c))} \sim (p(\ell_c) \cdot c)^{1/\beta} \rightarrow \lambda^{1/\beta} \text{ as } c \rightarrow \infty,
\end{align*} the second factor converges to \(\lambda^{1/\beta}\), and
it follows that \(Z \sim {\rm ML}(\beta,1)\). Since
\(c \rightarrow \infty\) is equivalent to \(\ell \uparrow x_R\), the
assertion follows with \(\ell:=\ell_c\) and \(p_\ell:=p(\ell_c)\). \qed

For a scale parameter \(\sigma > 0\), we write
\({\rm ML}(\beta, \sigma)\) for the distribution of \(\sigma Z_\beta\).
The Mittag-Leffler distribution with parameter \(\beta \in (0,1]\) is a
heavy-tailed positive distribution for \(\beta < 1\), with infinite
mean. However, as \(\beta \uparrow 1\), \({\rm ML}(\beta, \sigma)\)
converges weakly to the exponential distribution \({\rm Exp}(\sigma)\)
with mean \(\sigma\). This means that although its moments are all
infinite, the Mittag-Leffler distribution may (if \(\beta\) is close to
1) be indistinguishable from the exponential distribution, for the
purposes of applied statistics.

We caution the reader that, somewhat confusingly, there is another
distribution called the ``light-tailed'' Mittag-Leffler distribution.
This is in fact the limiting distribution of the renewal process
\(N(t)\) above (see Meerschaert and Scheffler (2004)). For a detailed
reference on the Mittag-Leffler distribution, see e.g.~ Haubold et al.
(2011), and for algorithms, see e.g.~the R package
\texttt{MittagLeffleR} (Gill and Straka, 2017).

\begin{description}
\item[\textbf{Remark:}]
If \(\beta = 1\), the result of the Theorem above is standard, see
Equation (2.2) in Gut and Hüsler (1999). In Anderson (1987) a similar
result is shown with a different choice of scaling constant. Meerschaert
and Stoev (2008) proved a limit theorem for the maxima of iid random
variables separated by infinite mean waiting times. They also show that
the hitting time of the limit process is Mittag-Leffler distributed.
Basrak and Špoljarić (2015) describe the asymptotic distribution, under
similar assumptions, of all upper orders statistics of the maximum
process using point processes.
\item[\textbf{Remark:}]
When \(0 < \beta < 1\), the renewal process \(N(t)\) is \emph{not
stationary}, and hence the results by Hsing et al. (1988) on the
exceedances of stationary sequences do not apply.
\end{description}

\hypertarget{statistical-inference}{%
\section{Statistical Inference}\label{statistical-inference}}

\hypertarget{mittag-leffler-stability-plots}{%
\subsection{Mittag-Leffler Stability
Plots}\label{mittag-leffler-stability-plots}}

We assume a dataset of the form \(\{(t_i, x_i), i \in I\}\), where
\(t_i\) are timestamps and \(x_i\) are magnitudes. In this article, we
focus on modelling the distribution of threshold \emph{exceedance times}
rather than the \emph{exceedances} themselves. For the latter, see the
POT approach (Beirlant et al., 2006, p. Ch.5.3; Coles, 2001, p. Ch.4.;
Davison and Smith, 1990; Embrechts et al., 2013 Ch. 6.5.; Leadbetter,
1991; Smith, 1984). With the previous section setting the stage, we
assume that only the large magnitudes (or ``shocks'') follow a MRP. The
smaller magnitudes are possibly nuisance data that may be irrelevant
with respect to the modelling of extremes and their waiting times.
Accordingly, we assume a threshold \(\ell\), and discard all data where
\(x_i \le \ell\), yielding the thresholded dataset
\(\{(t^{(\ell)}_i, x^{(\ell)}_i): i \in I(\ell)\}\), with
\(x^{(\ell)}_i=x_i-\ell\), where
\(I(\ell):=\left\{i \in I: x_i \ge \ell \right\} \subseteq I\).

However, this raises the question of how high or low the threshold
should be chosen. Recall that with the POT method for the modelling of
extremes, the distribution of the thresholded exceedances converges to a
GPD, and the threshold must be chosen high enough so that the
exceedances fit a GPD well. The threshold choice translates into a
bias-variance trade-off: low thresholds yield more data to fit a GPD
(low variance), but the distribution of exceedances may deviate from a
GPD distribution (high bias); and high thresholds present the opposite
scenario. It is understood, however, that there is a \emph{range} of
threshold values that ought to yield similar GPD parameter estimates. In
other words, plots of threshold vs.~parameter estimate ought to exhibit
that parameter estimates are \emph{stable} with respect to the choice of
threshold.

For inter-exceedance \emph{times} (IETs), the same idea applies. By the
Theorem, the inter-exceedance times converge to a Mittag-Leffler
distribution (MLD), and importantly, the tail parameter \(\beta\) is
\emph{independent} of (or \emph{stable} with respect to) the choice of
threshold. For high thresholds, the IET distribution may fit an MLD
well, but few IETs may be left for a fit (high variance). For low
thresholds, their distribution may be not well represented by an MLD;
or, too many events may be interpreted as shocks, resulting in a biased
estimate of the parameters of the limiting MLD. If we had MRP data
\(\{(W_i, J_i): i \in I\}\) as defined in Section \ref{sec:CTRE}, we
could apply any tail (and scale) estimator such as e.g.~Hill's estimator
(Hill, 1975) to the \(W_i\) and thus infer the distribution of
inter-exceedance times. However, in real-world datasets, the MRP
assumption is too strong for low thresholds. The unknown data generating
process is likely more complex, with dependencies or multiple data
generating mechanisms at low thresholds and structures that may vanish
for higher thresholds. Fitting only the IETs means that the iid
assumption applies to exceedances and IETs only.

\begin{figure}

{\centering \includegraphics[width=1\linewidth]{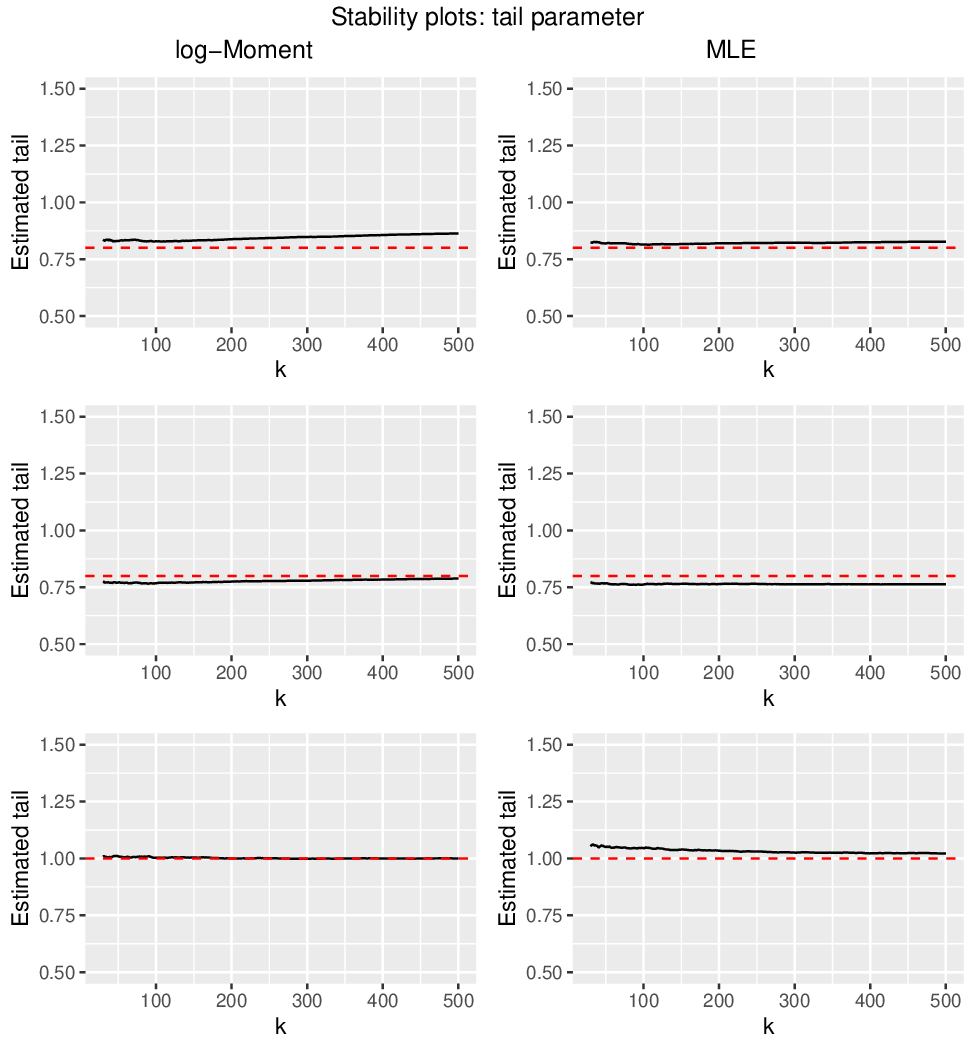} 

}

\caption{\label{Fig:TailSimu} Stability Plots for m=100 simulation runs for stable distributed waiting times with a tail parameter of 0.8 (top row), Pareto distributed waiting times with tail parameter 0.8 (middle row) and exponentially distributed waiting times (lower row). Left column: log-Moment estimator, right column: MLE. The grey thin lines are the stability plots for the different simulation runs and the dark lines are their means. The red dotted line shows the true tail parameter.}\label{fig:TailSimuplots}
\end{figure}

For high thresholds where few IETs are left, estimating a tail parameter
by fitting a MLD has advantages over e.g.~the Hill estimator: since the
underlying distribution of IETs \emph{is} a MLD, the fits are more
accurate. Figure \ref{fig:Hillplots} shows Hill plots for 200 IETs, for
two Mittag-Leffler datasets of sizes 200 and 10000 (and an irrelevant
distribution of magnitudes). Clearly the variance is too high to give
useful estimates. For corresponding fits based on the assumption of a
Mittag-Leffler distribution, see Figure \ref{Fig:TailSimu}, top row,
yielding more reliable predictions (more discussion further below).

\begin{figure}

{\centering \includegraphics[width=0.5\linewidth]{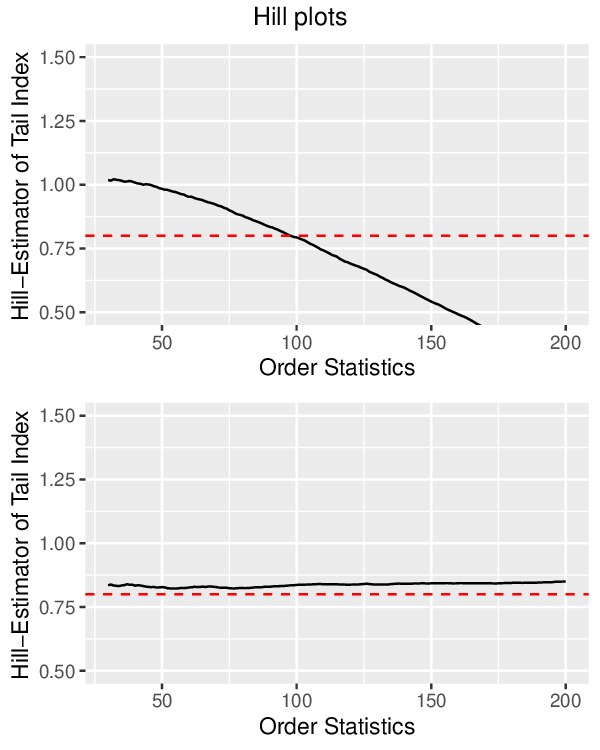} 

}

\caption{\label{fig:Hillplots} Hill plots for m=100 simulated Mittag-Leffler datasets with true tail 0.8 and sample size 200 (first panel) and 10000 (second panel), with number of upper order statistics r on which the Hill estimator is based on the x-axis. The grey thin lines are the Hill plots for the different simulation runs and the dark lines are their means. The red dotted line shows the true tail parameter. }\label{fig:Hillplots}
\end{figure}

\hypertarget{sec:ML}{%
\subsection{Fitting Mittag-Leffler Distributions}\label{sec:ML}}

Historically, the first method proposed for the estimation of the
Mittag-Leffler distribution parameters was the fractional moment
estimator by Kozubowski (2001). Unlike the first moments, the fractional
moments of order \(p\) for \(p<\beta\) exist and are tractable. One
drawback of this method is that constant priors for the tail parameter
are needed for the calculation of the estimates. Cahoy et al. (2010)
proposed a moment estimator of the log-transformed data, which does not
require any prior. Furthermore, they performed simulation studies
illustrating that the log-Moment outperforms the fractional moment
estimator with respect to bias and root mean squared error (RMSE).

The Maximum Likelihood Estimator (MLE) for the MLD is not
straightforward to implement, since the MLD density admits no closed
form analytical expression. In the R Package \texttt{MittagLeffleR}
(Gill and Straka, 2017), MLE is implemented via numerical optimization.
The MLE slightly outperforms the log-Moment estimator regarding bias and
RMSE for big enough sample sizes, but is computationally very intensive.
Figure \ref{fig:MSE} shows that both estimators perform well, even for
small sample sizes. For smaller tails both estimators show an increasing
RMSE for the scale estimation due to an increasing variance. This
results from the fact that in case of very small tails, single very
large values can occur.

\begin{figure}

{\centering \includegraphics[width=0.9\linewidth]{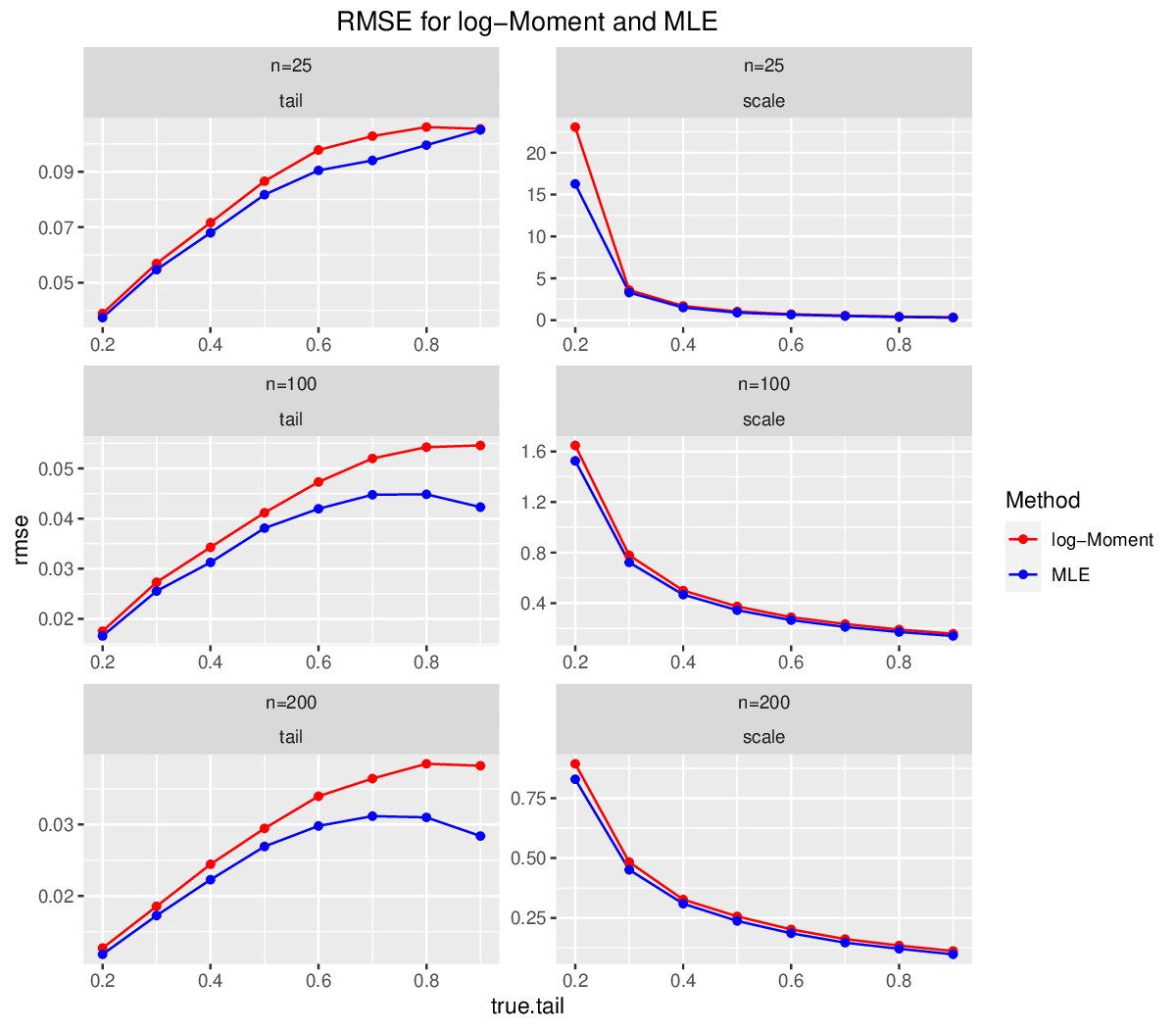} 

}

\caption{\label{fig:MSE} RMSE for the estimation of tail (left column) and scale (right column) parameters via log-Moment estimator and MLE of a Mittag-Leffler sample with varying sample size n=25, 100, 200, varying tails on the x-axis and fixed scale equal to one, based on 5000 simulation runs.}\label{fig:MSE}
\end{figure}

With the above comparisons of estimators for the MLD in mind, we propose
to use the log-Moment estimator, or the MLE when computational resources
are not an issue.

\hypertarget{weighing-the-evidence-for-non-exponential-inter-arrival-times}{%
\subsection{Weighing the evidence for non-exponential inter-arrival
times}\label{weighing-the-evidence-for-non-exponential-inter-arrival-times}}

\begin{figure}

{\centering \includegraphics[width=0.9\linewidth]{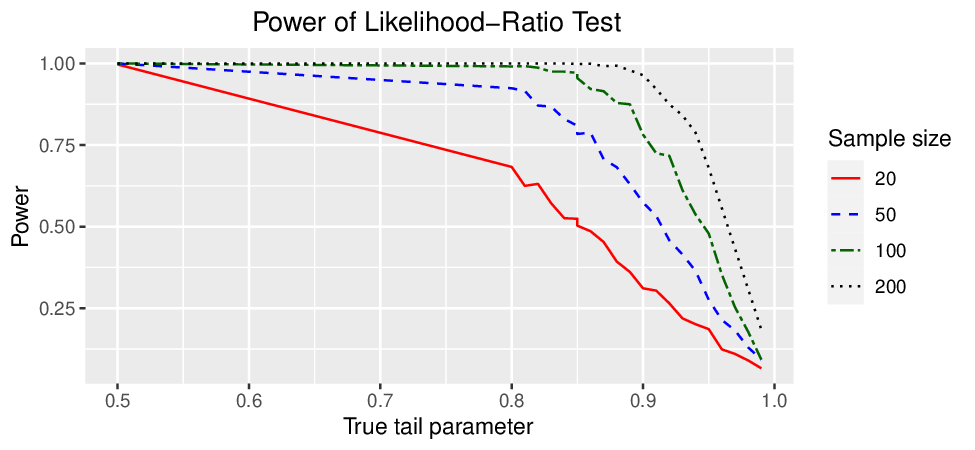} 

}

\caption{\label{Fig:LRT} Power for bootstrapped LRT for different sample sizes, varying tails and scale parameter equal to 1.}\label{fig:LRT_power}
\end{figure}

Since the exponential distribution is nested in the Mittag-Leffler
family of distributions, a Likelihood-ratio Test (LRT) seems to be an
appropriate way to choose between a model with exponential and
Mittag-Leffler inter-exceedance times. Although the two models are
nested, the asymptotic distribution is not \(\chi^2_1\)-distributed, and
Wilk's Theorem does not hold: under \(H_0\), the parameter \(\beta\) of
Mittag-Leffler distribution is equal to \(1\), and hence lies on the
boundary of the parameter space \((0,1]\). Instead, a valid approach is
a bootstrapped Likelihood-ratio test (see e.g.~Davison et al., 1997).
Figure \ref{Fig:LRT} displays the (simulated) power for the bootstrapped
LRT for Mittag-Leffler distributions with varying tail parameters based
on 1000 simulation runs. As expected, the power decreases for tail
parameters close to one, since the Mittag-Leffler distribution converges
as \(\beta \uparrow 1\) to an exponential distribution; it becomes hard
to differentiate a Mittag-Leffler distribution from an exponential.

\hypertarget{algorithm-for-the-inference-on-inter-exceedance-times}{%
\subsection{Algorithm for the inference on inter-exceedance
times}\label{algorithm-for-the-inference-on-inter-exceedance-times}}

The Theorem in Section \ref{sec:scaling} implies that for a high
threshold \(\ell\) we may approximate the distribution of \(T(\ell)\)
with an \({\rm ML}(\beta, b(1/p_{\ell}))\) distribution, where the
function \(b(c)\) varies regularly at \(\infty\) with parameter
\(1/\beta\). Building on the POT (Peaks-Over-Threshold) method, we
propose the following estimation procedure for the distribution of
inter-exceedance time \(T(\ell)\):

\begin{enumerate}
\def\labelenumi{\arabic{enumi}.}
\item
  Extract the \(K\) largest order statistics (i.e.~the \(K\) largest
  values, where e.g.~\(K = 300\)) \(x_{(1)}, \ldots, x_{(K)}\) together
  with their timestamps \(t_{(1)}, \ldots, t_{(K-1)}\).
\item
  Choose a minimum number of exceedances \(K_0\leq K\),
  e.g.~\(K_0 = 5\), and for each \(k\) ranging from \(K_0\) to \(K\):

  \begin{enumerate}
  \def\labelenumii{\alph{enumii})}
  \tightlist
  \item
    extract the set \(\mathcal T_k\) of exceedance times between the
    magnitudes exceeding the threshold \(X_{(k)}\)
  \item
    fit a Mittag-Leffler distribution to \(\mathcal T_k\), resulting in
    the parameter estimates \(\hat\beta_k\) and \(\hat \sigma_k\).
  \end{enumerate}
\item
  Plot \(k\) vs.~\(\hat \beta_k\). To the right (\(k \uparrow K\), low
  threshold), the asymptotics are off and bias is high. To the left
  (\(k \downarrow K_0\), high threshold), data is scarce and variance is
  high. In the middle, look for a region of \emph{stability} for a
  parameter estimate \(\hat \beta\). Choose as \(\hat \beta\) a
  representative value from this region.
\item
  Plot \(k\) vs.~\(k^{1/\hat \beta} \hat \sigma_k\). Again, choose a
  region of \emph{stability} and a representative \(\hat \sigma_0\) from
  that region.
\end{enumerate}

The inferred values \(\hat \beta\) and \(\hat \sigma_0\) can then be
interpreted as follows: Setting the threshold at \(X_{(k)}\), the
threshold exceedance times follow a Mittag-Leffler distribution with
shape parameter \(\hat \beta\) and scale parameter
\(\hat \sigma_0 k^{-1/\beta}\).

To clarify Step 4: Recall that by the theorem,
\(\sigma_k / b(1/p_{\ell})\) is expected to stabilize around a constant
as \(k\) decreases. Since \(b\) is regularly varying with parameter
\(1/\beta\), we have
\(b(1/p_{\ell}) = p_{\ell}^{-1/\beta} / L(1/p_{\ell})\) for some slowly
varying function \(L\). Approximating \(p_{\ell}\) by
\(\hat p_{\ell} := k / n\), we have
\[\text{const} \approx \sigma_k / b(1/p_{\ell}) = \sigma_k \times p_{\ell}^{1/\beta} L(1/p_{\ell}) \approx \sigma_k \times k^{1/\beta} n^{-1/\beta} L(n/k)\]
Assuming that the variation of \(L(n/k)\) is minor, we can hence see
that \(\sigma_k k^{1/\beta}\) stabilizes.

\begin{description}
\item[\textbf{Remark}:]
We approximated \(p_{\ell}\), the probability that an event is larger
than \(l\), by its relative frequency. One could also approximate this
tail probability via the GPD distribution fitted to the exceedances.
\end{description}

\begin{figure}
\includegraphics[width=\textwidth]{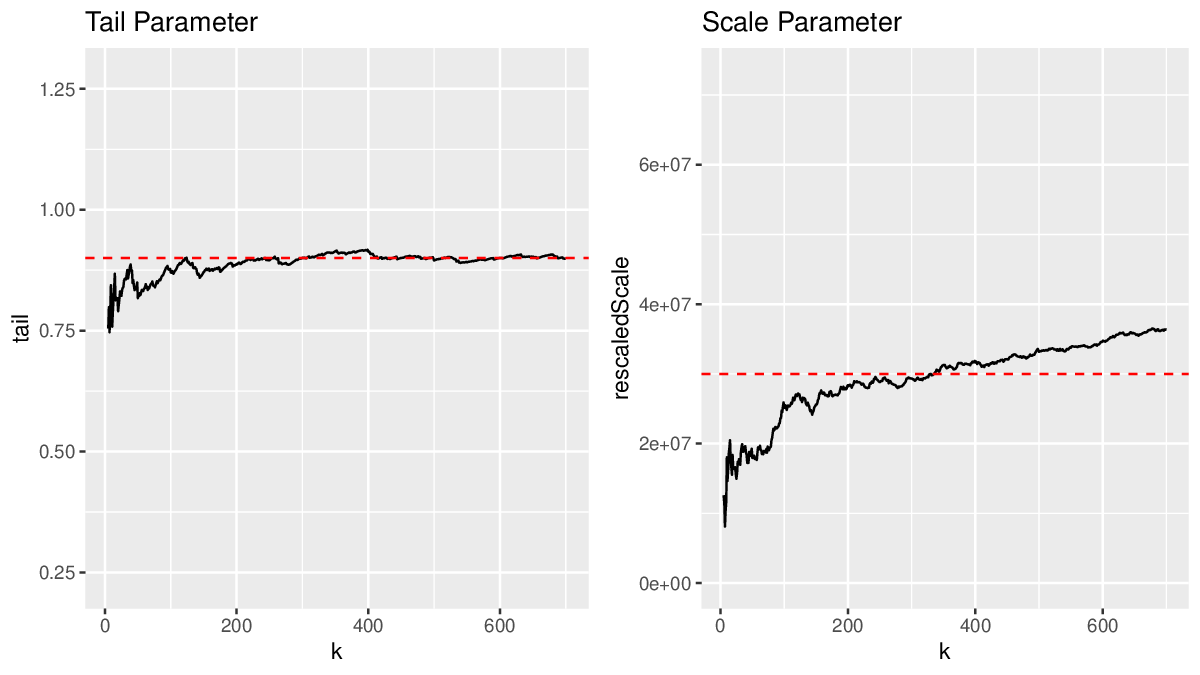} \caption{\label{fig:flares} Stability plots for the tail and scale parameter of the Mittag-Leffler distribution of the Solar Flare dataset. Dotted horizontal lines are at $\beta = 0.9$ and $\sigma_0 = 3 \times 10^7$ seconds $\approx 0.95$ years.}\label{fig:solar-flare-tail-scale}
\end{figure}

For computationally efficient estimates of the Mittag-Leffler parameters
we have used the method of log-transformed moments. This estimation
method provides point estimates as well as confidence intervals based on
sampling variance (Cahoy, 2013), and has been implemented in the R
software package \texttt{MittagLeffleR} (Gill and Straka, 2017). The
stability plots for \(\hat \beta\) and \(\hat \sigma_0\) can be
furnished with these confidence intervals, see e.g.~Figure
\ref{fig:flares}, to produce (non-simultaneous) confidence bands. These
stability plots were produced with the R package \texttt{CTRE} (Hees and
Straka, 2018). We have verified the validity of our estimation algorithm
via simulations, see Section \ref{Simulationstudy}.

\hypertarget{Simulationstudy}{%
\section{Simulation Study}\label{Simulationstudy}}

To test our inference method via stability plots, we have simulated
\(m=100\) datasets with \(n=10000\) independent waiting time and
magnitude pairs \((W_k, J_k)\) for waiting times that follow

\begin{enumerate}
\def\labelenumi{(\roman{enumi})}
\item
  a stable distribution,
\item
  a Pareto distribution and
\item
  an exponential distribution.
\end{enumerate}

The magnitudes are in all scenarios standard Gumbel distributed. In
order to have exact analytical values available for \(\beta\) and
\(\sigma_0\), a distribution for \(W_k\) needs to be chosen for which
\(b(n)\) from \eqref{eq:stability} is known.

\textbf{Case (i):} For (i) we choose \(W_k \stackrel{d}{=} D\), where
\(D\) is as in \eqref{eq:stability}, then due to the stability property
we have the \emph{equality} of distribution
\(W_1 + \ldots + W_n \stackrel{d}{=} b(n) D\), for
\(b(n) = n^{1/\beta}\). Using the parametrisation of Samorodnitsky and
Taqqu (1994), a few lines of calculation (see e.g.~the vignette on
parametrisation in Gill and Straka, 2017) show that \(D\) must have the
stable distribution \(S_\beta(\cos(\pi \beta/2)^{1/\beta}, +1, 0)\),
which is implemented in the R package \texttt{stabledist} by Wuertz et
al. (2016). By the Theorem, the distribution of \(T(\ell)\) is
approximately \[
{\rm ML}(\beta, p_{\ell}^{-1/\beta}) 
= {\rm ML}(\beta, k^{-1/\beta} n^{1/\beta}),
\] which means \(\sigma_0 = n^{1/\beta}\).

\textbf{Case (ii):} In the Pareto example we choose
\(P(W>t)=Ct^{-\beta}\) with \(C=(1/\Gamma(1-\beta))^{1/\beta}\). We have
chosen \(\beta=0.8\) in both cases (i) and (ii).

\textbf{Case (iii):} We choose exponentially distributed waiting times
with a rate parameter of \(1\).

Figure \ref{Fig:TailSimu} shows the graphical ``stability plots'' for
the estimation of the tail parameter, where

\begin{itemize}
\tightlist
\item
  rows correspond to cases (i), (ii) and (iii), and
\item
  columns correspond to estimators (log-Moment estimator,
  Maximum-Likelihood).
\end{itemize}

We plot the tail parameter estimates \(\hat \beta(k)\) against \(k\) for
each of the \(m=100\) simulation runs. Thin grey lines represent
individual simulation runs, and the thicker black line is their mean.
Recall that \(k\) is the index of the order statistics of \(J_k\) at
which the threshold \(\ell\) is placed.

\begin{figure}

{\centering \includegraphics[width=1\linewidth]{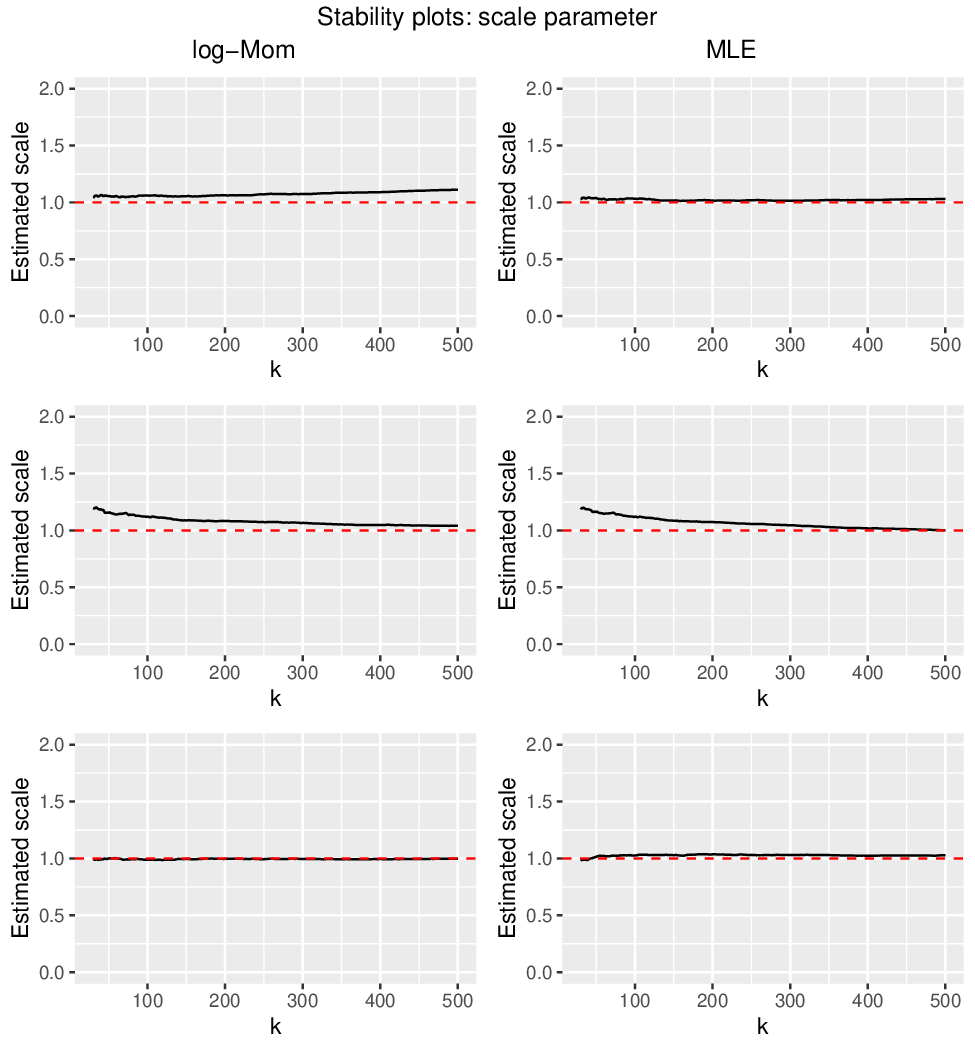} 

}

\caption{\label{Fig:ScaleSimu} Stability Plots for m=100 simulation runs for stable distributed waiting times with a tail parameter of 0.8 (top row), Pareto distributed waiting times with tail parameter 0.8 (middle row) and exponentially distributed waiting times (lower row). Left column: log-Moment estimator, right column: MLE. The grey thin lines are the stability plots for the different simulation runs and the dark lines are their means. The red dotted line shows the true scale parameter.}\label{fig:ScaleSimuplots}
\end{figure}

The performance of log-Moment and Maximum-Likelihood estimators for the
tail resp. the scale parameter is shown in Figure \ref{Fig:TailSimu}
resp. \ref{Fig:ScaleSimu}. Both estimators show good performance, with a
slight advantage for the Maximum-Likelihood estimator. This advantage is
paid for with a higher computational cost. Moreover, the bottom row in
Figure \ref{Fig:TailSimu} shows clearly that the log-Moment and
Maximum-Likelihood estimators generalize to the exponential case
(\(\beta = 1\)) whereas the Hill-estimator would not.\\
Notice that any data below a certain threshold is discarded in our
approach and hence need not satisfy the iid assumption. In our
simulation study we only simulated data from iid sequences of
\((W_k)_{k \in \mathbb{N}}\); in real data situations it is likely that
there are dependencies or more complex data generating processes which
vanish for higher thresholds. We only need the exceedances and IETs to
be iid.

\hypertarget{data-example}{%
\section{Data example}\label{data-example}}

We now want to apply the proposed method to a real data example, the
solar flare data which was already mentioned in Section 1 and can be
seen in Figure \ref{fig:thresholdedBursty}. The data were extracted from
the ``complete Hard X Ray Burst Spectrometer event list'', a
comprehensive reference for all measurements of the Hard X Ray Burst
Spectrometer on NASA's Solar Maximum Mission from the time of launch on
Feb 14, 1980 to the end of the mission in Dec 1989. 12,772 events were
detected, with the ``vast majority being solar flares''. To assure
stationarity and due to missing values during the years 1983 and 1984,
we based our analysis just on the year 1982, in which 2,488 events
happened. The list includes the start time, peak time, duration, and
peak rate of each event. We have used ``start time'' as the variable for
event times, and ``peak rate'' as the variable for event magnitudes.

Before we apply the approach described in Section 5 to the solar flare
data, we first have to check if all model assumptions are fulfilled. The
CTRE model is based on three main assumptions, which are repeated below.
For each assumption, we suggest one means of checking if it holds:

\begin{description}
\item[iid:]
After removing the ``noise observations'' below the smallest threshold
\(\ell_0\), the pair sequence \((T(\ell_0, i), X(\ell_0,i))\) is iid. An
indication if this is true is given by an auto-correlation plot. Since
we are expecting the inter-exceedance times to be Mittag-Leffler
distributed and hence to have infinite mean but finite log-moments, we
first take the logarithms of the times.
\item[Uncoupled:]
Each \(T(\ell, i)\) is independent of each \(X(\ell, i)\). We propose an
empirical copula plot to check for any dependence.
\item[\({\rm ML}(\beta, \sigma)\) distribution of \(T(\ell, i)\):]
Apply a cutoff at the lowest threshold \(\ell_0\), extract the threshold
crossing times, and create a QQ Plot for the Mittag-Leffler
distribution. Use a Log-moment estimate of the tail parameter for the
theoretical / population quantiles of the plot.
\item[\textbf{Remark:}]
The ACF plots of course can just give an indication whether there are
dependencies, since they actually just measure linear dependencies.
Furthermore, if one calculates the ACF for the logarithmic
inter-exceedance times, the ACF indicates on the original scale a
multiplicative dependence.
\end{description}

Figures \ref{fig:flare-diagnostics-1}, \ref{fig:flare-diagnostics-2} and
\ref{fig:flare-diagnostics-3} show the diagnostic plots for a minimum
threshold chosen at the 100th order statistic. There is some residual
autocorrelation for the sequence of threshold exceedance times that is
not accounted for by the CTRE model.

\begin{figure}
\includegraphics[width=\textwidth]{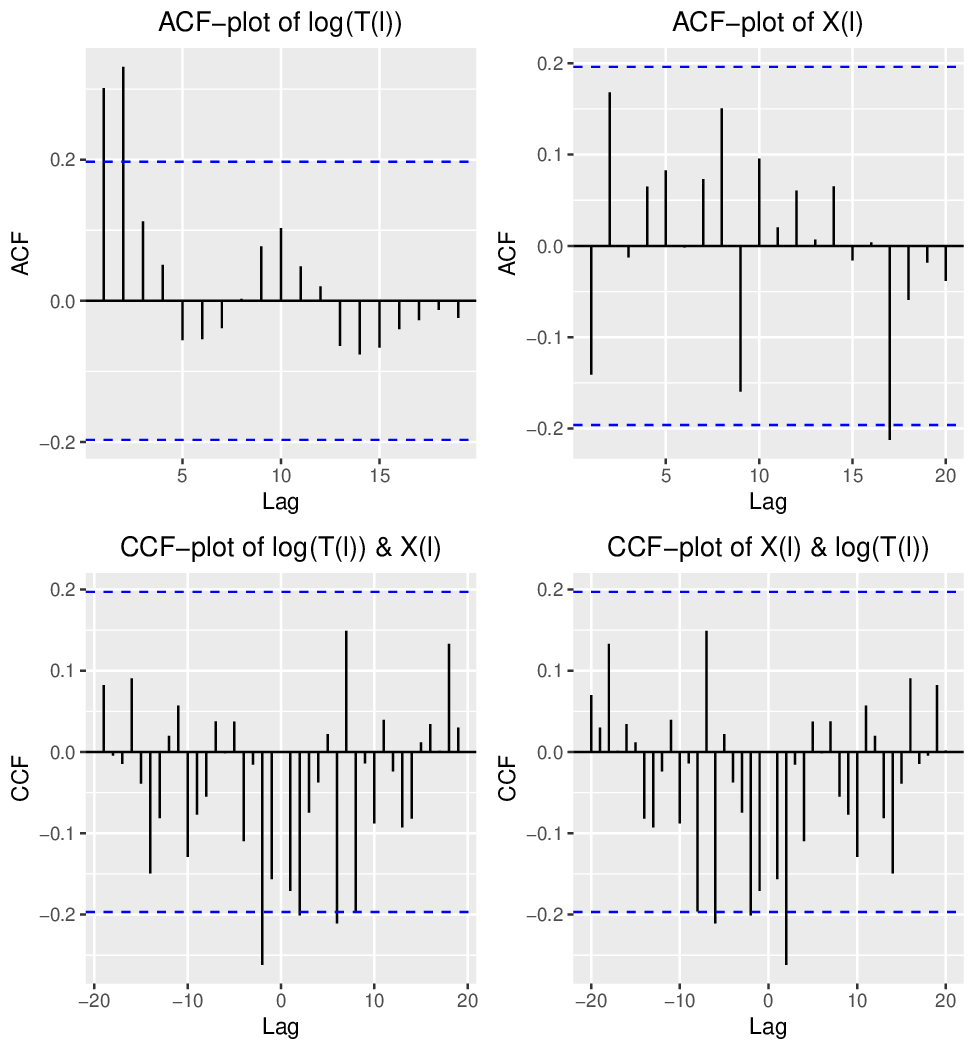} \caption{\label{fig:flare-diagnostics-1} Diagnostic plots for the solar flare data based on the 100 upper order statistics: auto-correlation and cross-correalation function.}\label{fig:flare-diagnostics-1}
\end{figure}

\begin{figure}

{\centering \includegraphics[width=0.7\linewidth]{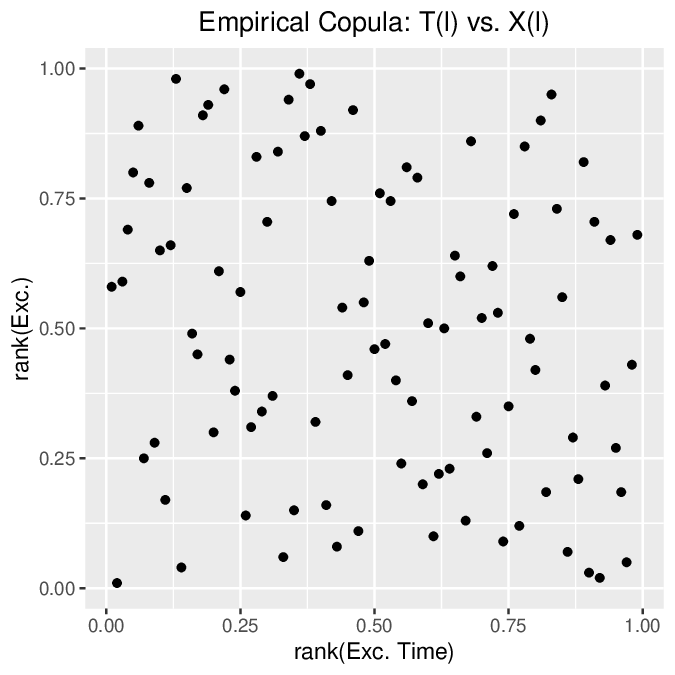} 

}

\caption{\label{fig:flare-diagnostics-2} Diagnostic plots for the solar flare data: empirical copula.}\label{fig:flare-diagnostics-2}
\end{figure}

\begin{figure}

{\centering \includegraphics[width=0.7\linewidth]{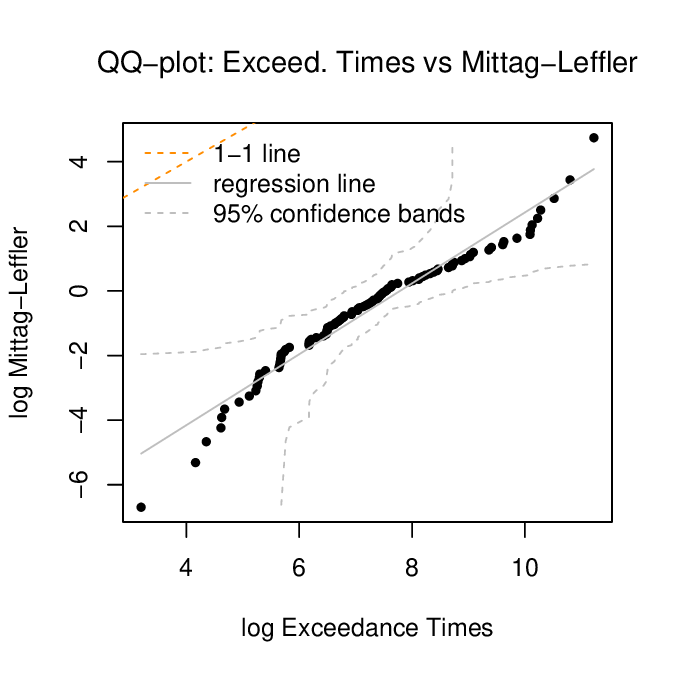} 

}

\caption{\label{fig:flare-diagnostics-3} Diagnostic plots for the solar flare data: QQ Plot.}\label{fig:flare-diagnostics-3}
\end{figure}

Figure \ref{fig:flares} shows the stability plots for the solar flare
data, on the left for the tail parameter and on the right for the scale
parameter. The dark grey ranges correspond to 95\% confidence intervals,
which are derived from the asymptotic normality of the Log-moment
estimators (Cahoy, 2013) and the \(\delta\)-method (Gill and Straka,
2017); dashed lines show the deduced true values of \(\beta\)
resp.~\(\sigma_0\). The stability plot for the tail stabilizes nicely
around 0.9 (dashed line), while the scale parameter stabilizes less
obviously near \(3 \times 10^7\) (dashed line). The growth of the scale
parameter for lower threshold appears to be closer to linear in
\(p_{\ell}\), rather than proportional to \(p_{\ell}^{1/0.9}\) as
suggested by the Mittag-Leffler fits. The reason for this is likely that
the overall goodness of fit as compared to an exponential distribution
is improved due to the peaked shape of the Mittag-Leffler distribution
near \(0\), rather than its tail behaviour at \(\infty\). The reported
fit should hence come with the caveat that a Mittag-Leffler distribution
models exceedance times well only up to certain time-scales. More
research is needed into the modelling of scale transitions, where
inter-exceedance times appear to have different power laws across
different time scales.

The fit with a Mittag-Leffler distribution (\(\beta = 0.9\)) seems to be
good (see Figure \ref{fig:flare-diagnostics-3}), although there are
signs that the power-law tail tapers off for very large inter-threshold
crossing times. There is no apparent dependence between threshold
exceedance times and event magnitudes seen in the copula plot (see
Figure \ref{fig:flare-diagnostics-2}). We also conduct a bootstrapped
LRT for the null hypothesis of exponentially distributed inter-arrival
times and received a \(p\)-value of \(p<0.01\).

\hypertarget{predicting-the-time-of-the-next-threshold-crossing}{%
\section{Predicting the time of the next threshold
crossing}\label{predicting-the-time-of-the-next-threshold-crossing}}

According to Figure \ref{fig:flares}, for a threshold \(\ell\) at the
\(k\)-th order statistic, the estimated threshold exceedance time
distribution is approximately \[
T(\ell) \sim {\rm ML}(\hat{\beta}, k^{-1/\hat{\beta}} \hat{\sigma}_0), 
\] where \(\hat{\beta} = 0.9\) and
\(\hat{\sigma}_0 = 3.0 \times 10^7 {\rm sec}\). Unlike the exponential
distribution, the Mittag-Leffler distribution is not memoryless, and the
probability density of the time \(t\) until the next threshold crossing
will depend on the time \(t_0\) elapsed since the last threshold
crossing. This density is approximately equal to \[
p(t|\beta, \sigma_0, \ell, t_0) = \frac{f(t + t_0 | \beta, k^{-1/\beta} \sigma_0)}{\mathbf P[T_\ell > t_0]}
\] where \(f(\,\cdot\, | \beta, k^{-1/\beta} \sigma_0)\) is the
probability density of \({\rm ML}(\beta, k^{-1/\beta} \sigma_0)\). The
more time has passed without a threshold crossing, the more the
probability distribution shifts towards larger values for the next
crossing (see Figure \ref{fig:hazard}, left panel). The hazard rate \[
h(t) = \frac{f(t| \beta, k^{-1/\beta} \sigma_0))}{\int_t^\infty f(\tau| \beta, k^{-1/\beta} \sigma_0))\,d\tau}
\] approximates the risk of a threshold crossing per unit time, and is a
decreasing function for the Mittag-Leffler distribution.

\begin{figure}
\includegraphics[width=\textwidth]{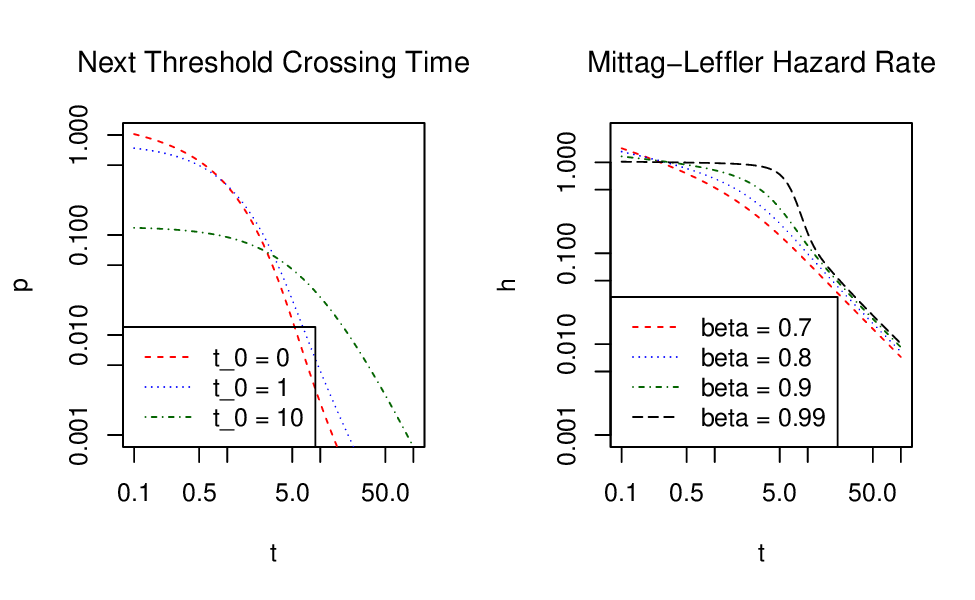} \caption{\label{fig:hazard} Left: Conditional distribution of time until next threshold crossing, depending on elapsed time $t_0$ since last crossing ($\beta = 0.9$, $\sigma_0 = 1$). Right: Hazard rate depending on tail parameter $\beta$.}\label{fig:hazard}
\end{figure}

The closer \(\beta\) is to \(1\), the more the hazard rate mimics that
of an exponential distribution (a constant function, see Figure
\ref{fig:hazard}, right panel).

\hypertarget{discussion-conclusion}{%
\section{Discussion \& Conclusion}\label{discussion-conclusion}}

We proposed a new model and inference procedure for the inter-exceedance
times of ``bursty'' time series, which have been studied intensively in
statistical physics. Burstiness is characterized by power-law waiting
times between events, and we have shown that the Mittag-Leffler
distribution arises naturally as a scaling limit for the
inter-exceedance times of high thresholds. Moreover, we have derived the
following non-linear scaling behaviour:
\(\sigma \sim p_{\ell}^{-1/\beta}\), where \(\sigma\) is the scale
parameter of the distribution of threshold exceedance times,
\(p_{\ell}\) is the fraction of magnitudes above the threshold \(\ell\),
and \(\beta\) the exponent of the power law. This ``anomalous'' scaling
behaviour in the bursty setting entails two phenomena:

\begin{enumerate}
\def\labelenumi{\roman{enumi})}
\item
  a heavy tail of the inter-arrival time distribution of threshold
  crossings (long rests), and
\item
  a high propensity for more threshold crossing events immediately after
  each threshold crossing event (bursts).
\end{enumerate}

The Mittag-Leffler distribution captures both phenomena, due to its
heavy tail as well as its stretched exponential (peaked) asymptotics for
small times. It generalizes the exponential distribution, and in the
solar flare data example, this generalization is warranted, because the
likelihood-ratio test is strongly significant.

When we introduced the CTRE model, we assumed that all events are iid.
This assumption is likely sufficient but not necessary for our limit
theorem to hold. Moreover, any data below a (minimum) threshold
\(\ell_0\) is discarded for CTREs, and hence need not satisfy the iid
assumption. For the purposes of statistical inference, we merely require
that the IETs are iid.

The CTRE approach to model ``non-Poissonian'' threshold crossing times
should be contrasted with the well-documented approach of clusters of
extremes, see e.g.~Ferro and Segers (2003). When the underlying
stochastic process is \emph{stationary}, the exceedances of high
thresholds form, asymptotically, a \emph{Cluster Poisson Process}. This
result was established in (Hsing et al., 1988). In the setting of this
article, however, clustering-like dynamics occur due to the
\emph{non-stationarity} of the underlying renewal-reward process, which
has infinite-mean renewal times. Hence CTRE and Cluster Poisson Process
should not be viewed as competing methods, as the underlying data
generating processes are quite different. To differentiate between the
two models we use the term \emph{bursty} which is standard in the
context of heavy-tailed inter-arrival times in the physics community
(e.g.~Barabási, 2005; Karsai et al., 2012; Vajna et al., 2013; Vasquez
et al., 2006)). In weighing the evidence for either of the two data
generating processes, criteria could be developed, based e.g.~on
measures of surprise (Lee et al., 2015), which may prove to be valuable
for future applied statisticians.

Finally, we note that assuming a purely scale-free pattern for event
times may be too rigid an assumption, which unnecessarily limits the
applicability of CTREs. Often, the heavy-tailed character of the
inter-arrival time distribution holds at short to intermediate time
scales, and is truncated (or tempered, reverting to an exponential
distribution) at very long time scales (see e.g.~Meerschaert et al.,
2012; and Aban et al., 2006). In such situations, a ``tempered''
Mittag-Leffler distribution may provide a better fit, which we aim to
introduce in follow-up work.

\hypertarget{acknowledgements}{%
\section*{Acknowledgements}\label{acknowledgements}}
\addcontentsline{toc}{section}{Acknowledgements}

Peter Straka was supported by the Discovery Early Career Research Award
DE160101147 on the Project ``Predicting Extremes when Events Occurin
Bursts'' by the Australian Research Council. Katharina Hees was
supported by the DAAD co-financed by the German Federal Ministry of
Education and Research (BMBF). The authors would like to thank
Prof.~Peter Scheffler for insights on stochastic process limits for
CTRMs, Prof.~Roland Fried for discussion regarding the statistical
methods and Gurtek Gill who helped create the MittagLeffleR R-package.

\newpage

\hypertarget{references}{%
\section*{References}\label{references}}
\addcontentsline{toc}{section}{References}

\hypertarget{refs}{}
\leavevmode\hypertarget{ref-Aban06}{}%
Aban, I.B., Meerschaert, M.M., Panorska, A.K., 2006. Parameter
estimation for the truncated pareto distribution. J. Am. Stat. Assoc.
101, 270--277.

\leavevmode\hypertarget{ref-Anderson1987}{}%
Anderson, K.K., 1987. Limit Theorems for General Shock Models with
Infinite Mean Intershock Times. J. Appl. Probab. 24, 449--456.

\leavevmode\hypertarget{ref-Bagrow2013}{}%
Bagrow, J.P., Brockmann, D., 2013. Natural emergence of clusters and
bursts in network evolution. Phys. Rev. X 3, 1--6.

\leavevmode\hypertarget{ref-Barabasi2005}{}%
Barabási, A.L., 2005. The origin of bursts and heavy tails in human
dynamics. Nature 435, 207--211.

\leavevmode\hypertarget{ref-Basrak2014}{}%
Basrak, B., Špoljarić, D., 2015. Extremes of random variables observed
in renewal times. Stat. Probab. Lett. 97, 216--221.

\leavevmode\hypertarget{ref-beirlantBook}{}%
Beirlant, J., Goegebeur, Y., Segers, J., Teugels, J., 2006. Statistics
of extremes: theory and applications. John Wiley \& Sons.

\leavevmode\hypertarget{ref-Benson2007}{}%
Benson, D.A., Schumer, R., Meerschaert, M.M., 2007. Recurrence of
extreme events with power-law interarrival times. Geophys. Res. Lett.
34.

\leavevmode\hypertarget{ref-Cahoy2013}{}%
Cahoy, D.O., 2013. Estimation of Mittag-Leffler Parameters. Commun.
Stat. - Simul. Comput. 42, 303--315.

\leavevmode\hypertarget{ref-Cahoy2010}{}%
Cahoy, D.O., Uchaikin, V.V., Woyczynski, W.A., 2010. Parameter
estimation for fractional Poisson processes. J. Stat. Plan. Inference
140, 3106--3120.

\leavevmode\hypertarget{ref-ColesBook}{}%
Coles, S., 2001. An Introduction to Statistical Modelling of Extreme
Values. Springer-Verlag, London.

\leavevmode\hypertarget{ref-davison1997bootstrap}{}%
Davison, A.C., Hinkley, D.V., others, 1997. Bootstrap methods and their
application. Cambridge university press.

\leavevmode\hypertarget{ref-davison1990models}{}%
Davison, A.C., Smith, R.L., 1990. Models for exceedances over high
thresholds. Journal of the Royal Statistical Society: Series B
(Methodological) 52, 393--425.

\leavevmode\hypertarget{ref-HXRBS}{}%
Dennis, B.R., Orwig, L.E., Kennard, G.S., Labow, G.J., Schwartz, R.A.,
Shaver, A.R., Tolbert, A.K., 1991. The complete Hard X Ray Burst
Spectrometer event list, 1980-1989.

\leavevmode\hypertarget{ref-embrechts2013modelling}{}%
Embrechts, P., Klüppelberg, C., Mikosch, T., 2013. Modelling extremal
events: For insurance and finance. Springer Science \& Business Media.

\leavevmode\hypertarget{ref-Esary1973}{}%
Esary, J.D., Marshall, A.W., 1973. Shock Models and Wear Processes 1,
627--649.

\leavevmode\hypertarget{ref-ferro2003inference}{}%
Ferro, C.A., Segers, J., 2003. Inference for clusters of extreme values.
Journal of the Royal Statistical Society: Series B (Statistical
Methodology) 65, 545--556.

\leavevmode\hypertarget{ref-MittagLeffleR}{}%
Gill, G., Straka, P., 2017. MittagLeffleR: Using the mittag-leffler
distributions in r.

\leavevmode\hypertarget{ref-gnedenko1983limit}{}%
Gnedenko, B., 1983. On limit theorems for a random number of random
variables, in: Probability Theory and Mathematical Statistics. Springer,
pp. 167--176.

\leavevmode\hypertarget{ref-Gut1999}{}%
Gut, A., Hüsler, J., 1999. Extreme Shock Models. Extremes 2, 295--307.

\leavevmode\hypertarget{ref-Haubold11}{}%
Haubold, H.J., Mathai, A.M., Saxena, R.K., 2011. Mittag-Leffler
Functions and Their Applications. J. Appl. Math. 2011, 1--51.

\leavevmode\hypertarget{ref-hawkes1971point}{}%
Hawkes, A.G., 1971. Point spectra of some mutually exciting point
processes. J. R. Stat. Soc. Ser. B Stat. Methodol. 438--443.

\leavevmode\hypertarget{ref-hees2016joint}{}%
Hees, K., Scheffler, H.-P., 2018a. On joint sum/max stability and
sum/max domains of attraction. Probability and Mathematical Statistics
38.

\leavevmode\hypertarget{ref-hees2017coupled}{}%
Hees, K., Scheffler, H.-P., 2018b. Coupled continuous time random
maxima. Extremes 21, 235--259.

\leavevmode\hypertarget{ref-CTRE}{}%
Hees, K., Straka, P., 2018. CTRE: Thresholding bursty time series.

\leavevmode\hypertarget{ref-hill1975simple}{}%
Hill, B.M., 1975. A simple general approach to inference about the tail
of a distribution. The annals of statistics 1163--1174.

\leavevmode\hypertarget{ref-Hsing88}{}%
Hsing, T., Hüsler, J., Leadbetter, M.R., 1988. On the exceedance point
process for a stationary sequence. Probability theory and related fields
78, 97--112.

\leavevmode\hypertarget{ref-Karsai2012}{}%
Karsai, M., Kaski, K., Barabási, A.L., Kertész, J., 2012. Universal
features of correlated bursty behaviour. Sci. Rep. 2.

\leavevmode\hypertarget{ref-Karsai2011}{}%
Karsai, M., Kivelä, M., Pan, R.K., Kaski, K., Kertész, J., Barabási,
A.L., Saramäki, J., 2011. Small but slow world: How network topology and
burstiness slow down spreading. Phys. Rev. E - Stat. Nonlinear, Soft
Matter Phys. 83, 1--4.

\leavevmode\hypertarget{ref-kozubowski2001}{}%
Kozubowski, T.J., 2001. Fractional moment estimation of linnik and
mittag-leffler parameters. Mathematical and computer modelling 34,
1023--1035.

\leavevmode\hypertarget{ref-Laskin2003}{}%
Laskin, N., 2003. Fractional Poisson process. Commun. Nonlinear Sci.
Numer. Simul. 8, 201--213.

\leavevmode\hypertarget{ref-leadbetter1991basis}{}%
Leadbetter, M.R., 1991. On a basis for `peaks over threshold'modeling.
Statistics \& Probability Letters 12, 357--362.

\leavevmode\hypertarget{ref-Lee15}{}%
Lee, J., Fan, Y., Sisson, S.A., 2015. Bayesian threshold selection for
extremal models using measures of surprise. Comput. Stat. Data Anal. 85,
84--99.

\leavevmode\hypertarget{ref-Meerschaert2010b}{}%
Meerschaert, M.M., Nane, E., Vellaisamy, P., 2011. The fractional
Poisson process and the inverse stable subordinator. Electron. J.
Probab. 16, 1600--1620.

\leavevmode\hypertarget{ref-MeerschaertRoyQin}{}%
Meerschaert, M.M., Roy, P., Shao, Q., 2012. Parameter estimation for
exponentially tempered power law distributions. Commun. Stat. - Theory
Methods 41, 1839--1856.

\leavevmode\hypertarget{ref-limitCTRW}{}%
Meerschaert, M.M., Scheffler, H.-P., 2004. Limit Theorems for
Continuous-Time Random Walks with Infinite Mean Waiting Times. J. Appl.
Probab. 41, 623--638.

\leavevmode\hypertarget{ref-MeerschaertSikorskii}{}%
Meerschaert, M.M., Sikorskii, A., 2012. Stochastic models for fractional
calculus. Walter de Gruyter.

\leavevmode\hypertarget{ref-MeerschaertStoev08}{}%
Meerschaert, M.M., Stoev, S.A., 2008. Extremal limit theorems for
observations separated by random power law waiting times. J. Stat. Plan.
Inference 139, 2175--2188.

\leavevmode\hypertarget{ref-Min2010}{}%
Min, B., Goh, K.I., Vazquez, A., 2011. Spreading dynamics following
bursty human activity patterns. Phys. Rev. E - Stat. Nonlinear, Soft
Matter Phys. 83, 2--5.

\leavevmode\hypertarget{ref-Oliveira2005}{}%
Oliveira, J., Barabási, A.L., 2005. Darwin and Einstein correspondence
patterns. Nature 437, 1251.

\leavevmode\hypertarget{ref-Omi2011}{}%
Omi, T., Shinomoto, S., 2011. Optimizing Time Histograms for
Non-Poissonian Spike Trains. Neural Comput. 23, 3125--3144.

\leavevmode\hypertarget{ref-R}{}%
R Core Team, 2018. R: A language and environment for statistical
computing. R Foundation for Statistical Computing, Vienna, Austria.

\leavevmode\hypertarget{ref-SamorodnitskyTaqqu}{}%
Samorodnitsky, G., Taqqu, M.S., 1994. Stable Non-Gaussian Random
Processes: Stochastic Models with Infinite Variance, Stochastic
modeling. Chapman Hall, London.

\leavevmode\hypertarget{ref-Sumita1983}{}%
Shanthikumar, J.G., Sumita, U., 1983. General shock models associated
with correlated renewal sequences. J. Appl. Probab. 20, 600--614.

\leavevmode\hypertarget{ref-Sumita1984}{}%
Shanthikumar, J.G., Sumita, U., 1984. Distribution Properties of the
System Failure Time in a General Shock Model. Adv. Appl. Probab. 16,
363--377.

\leavevmode\hypertarget{ref-Sumita1985}{}%
Shanthikumar, J.G., Sumita, U., 1985. A class of correlated cumulative
shock models. Adv. Appl. Probab. 17, 347--366.

\leavevmode\hypertarget{ref-Silvestrov2002a}{}%
Silvestrov, D.S., 2002. Limit Theorems for Randomly Stopped Stochastic
Processes. Springer (Berlin, Heidelberg).

\leavevmode\hypertarget{ref-ST04}{}%
Silvestrov, D.S., Teugels, J.L., 2004. Limit theorems for mixed max-sum
processes with renewal stopping. Ann. Appl. Probab. 14, 1838--1868.

\leavevmode\hypertarget{ref-smith1984threshold}{}%
Smith, R.L., 1984. Threshold methods for sample extremes, in:
Statistical Extremes and Applications. Springer, pp. 621--638.

\leavevmode\hypertarget{ref-stindl2018likelihood}{}%
Stindl, T., Chen, F., 2018. Likelihood based inference for the
multivariate renewal hawkes process. Computational Statistics \& Data
Analysis 123, 131--145.

\leavevmode\hypertarget{ref-Vajna2013}{}%
Vajna, S., Tóth, B., Kertész, J., 2013. Modelling bursty time series.
New J. Phys. 15, 103023.

\leavevmode\hypertarget{ref-Vasquez2006}{}%
Vasquez, a, Oliveira, J.G., Dezso, Z., Goh, K.-I., Kondor, I., Barabási,
A.L., 2006. Modeling bursts and heavy tails in human dynamics. Phys.
Rev. E 73, 361271--3612718.

\leavevmode\hypertarget{ref-Vazquez2007}{}%
Vazquez, A., Rácz, B., Lukács, A., Barabási, A.L., 2007. Impact of
non-poissonian activity patterns on spreading processes. Phys. Rev.
Lett. 98, 1--4.

\leavevmode\hypertarget{ref-wheatley2016hawkes}{}%
Wheatley, S., Filimonov, V., Sornette, D., 2016. The hawkes process with
renewal immigration \& its estimation with an em algorithm.
Computational Statistics \& Data Analysis 94, 120--135.

\leavevmode\hypertarget{ref-stabledist}{}%
Wuertz, D., Maechler, M., members., R. core team, 2016. Stabledist:
Stable distribution functions.

\end{document}